\documentclass{article}%
\usepackage{amssymb}
\usepackage{graphicx}
\usepackage{amsmath}
\usepackage{amsfonts}%
\providecommand{\U}[1]{\protect\rule{.1in}{.1in}}
\newtheorem{theorem}{Theorem}

\newtheorem{corollary}[theorem]{Corollary}

\newtheorem{proposition}[theorem]{Proposition}
\newtheorem{remark}[theorem]{Remark}

\begin{document}

\begin{center}
{\Large \textbf{The limiting spectral distribution in terms of spectral
density}}

Costel Peligrad and Magda Peligrad
\end{center}

\bigskip

(Random Matrices Theory and Applications vol 5 Nr 1 2016)

\bigskip

Department of Mathematical Sciences, University of Cincinnati, PO Box 210025,
Cincinnati, OH 45221-0025, USA.

Email: peligrc@ucmail.uc.edu; peligrm@ucmail.uc.edu

\bigskip

Key words: random matrices, correlated entries, limiting spectral
distribution, stationary random fields, spectral density.

Mathematics Subject Classification (2010): 60F15, 60G60, 60G10, 62E20.
Secondary 42A38

\begin{center}
\bigskip

\bigskip\textbf{Abstract}
\end{center}

For a large class of symmetric random matrices with correlated entries,
selected from stationary random fields of centered and square integrable
variables, we show that the limiting distribution of eigenvalue counting
measure always exists and we describe it via an equation satisfied by its
Stieltjes transform. No rate of convergence to zero of correlations is
imposed, therefore the process is allowed to have long memory. In particular,
if the symmetrized matrices are constructed from stationary Gaussian random
fields which have spectral density, the result of this paper gives a complete
solution to the limiting eigenvalue distribution. More generally, for matrices
whose entries are functions of independent identically distributed random
variables the result also holds.

\section{Introduction}

This paper is a contribution to the limiting spectral distribution for
symmetric matrices with correlated entries. Among the first results in this
direction are papers by Khorunzhy and Pastur \cite{KP}, Boutet de Monvel and
Khorunzhy \cite{BK}, who treated Gaussian random fields with absolutely
summable covariances. Khorunzhy \cite{K} considered matrices with correlated
entries imposing rates of convergence on some mixing coefficients, without
assuming the variables are Gaussian. On the other hand, there is interest in
studying linear filters of independent random variables as entries of a
matrix. Anderson and Zeitouni \cite{AZ} considered symmetric matrices with
entries that are linear processes of finite range of independent random
variables. In all the papers mentioned above the correlation between variables
are diminishing with time at certain polynomial rates. Such a dependence is
considered of weak type, since distant variables have weak interactions.

Very recently, results in several papers indicate that the weak dependence is
not needed for the existence of the limiting spectral distribution, and much
milder regularity conditions, or ergodic type, are responsible for this line
of problems. More precisely, Chakrabarty \textit{et al. }\cite{chak}
considered symmetrized random matrices selected from a large class of
stationary Gaussian random fields and argued that the limiting spectral
density always exists.

On the other hand, in two recent papers Banna \textit{et al.} \cite{BMP}, and
Merlev\`{e}de and Peligrad \cite{MP} showed that this type of general result
is not restricted to Gaussian fields. In \cite{BMP} was studied symmetric
random matrices whose entries are functions of zero mean square integrable
independent and identically distributed (i.i.d.) real-valued random variables.
Such kind of processes provide a general framework for stationary and ergodic
random fields. For this case, by using the substitution method, in \cite{BMP},
the study of the empirical eigenvalue distribution was reduced to the same
problem for a matrix with Gaussian entries and the same covariance structure,
without any other additional assumption. This universality result was combined
with a known result for the Gaussian case in \cite{KP} to extend it to
stationary random fields. However, because the result for the Gaussian case
given in \cite{KP} applies only to some weakly dependent random fields, the
characterization of the limit in \cite{KP} is restricted to the case when
covariances are absolutely summable.

In this paper we obtain a characterization of the limiting empirical spectral
distribution for symmetric matrices with entries selected from a stationary
Gaussian field under the sole condition that its spectral density exists. This
result, combined with the universality result in Banna \textit{et al.}
\cite{BMP} shows that for random matrices with entries functions of i.i.d.,
the limiting empirical spectral distribution exists and is characterized via
an equation satisfied by its Stieltjes transform which involves the spectral
density of the field. Applications of this result to linear and nonlinear
filters of a stationary random field are pointed out.

A general characterization of the limiting spectral distribution in terms of
the spectral density is also expected to hold for the covariance (Gram)
matrices whose entries are selected from stationary random fields. This fact
is suggested by results in papers by Boutet de Monvel \textit{et al.}
\cite{BKV}, Hachem \textit{et al. }\cite{HLN}, Bai and Zhou \cite{BZ}, Yao
\cite{Yao}, Banna \textit{et al.} \cite{BMP}, Merlev\`{e}de and Peligrad
\cite{MP}, among many others. The study of covariance matrices is beyond the
scope of this paper.

Our paper is organized as follows. In Section 2, we give definitions and state
the main results. Section 3 contains the proofs. Section 4 is dedicated to examples.

\section{Results}

Here are some notations used throughout the paper. For a matrix $A$, we denote
by $\mathrm{Tr}(A)$ its trace. We shall use the notation $\Vert X\Vert_{r}$
for the ${\mathbb{L}}_{r}$-norm ($r\geq1$) of a real valued random variable
$X,$ namely $\Vert X\Vert_{r}^{r}=\mathbb{E(}|X|^{r})$. For a set $B$ we
denote by $B^{\prime}$ its complement. For the convergence in distribution we
use the notation $\Rightarrow.$ The Lebesgue measure on $\mathbb{R}$ will be
denoted by $\lambda.$ The set of complex numbers with positive imaginary part
is denoted by $\mathbb{C}^{+}.$

For any square matrix $A_{n}$ of order $n$ with real eigenvalues $\lambda
_{1}\leq\dots\leq\lambda_{n}$, its empirical spectral measure and its
empirical spectral distribution function are respectively defined by
\[
\nu_{A_{n}}=\frac{1}{n}\sum_{k=1}^{n}\delta_{\lambda_{k}}\ \text{ and }%
\ F_{n}^{A_{n}}(x)=\frac{1}{n}\sum_{k=1}^{n}\mathbf{1}_{\{\lambda_{k}\leq
x\}}\,.
\]
The Stieltjes transform of $F_{n}^{A_{n}}$ is given by
\[
S_{A_{n}}(z)=\int\frac{1}{x-z}\mathrm{d}F_{n}^{A_{n}}(x)=\frac{1}%
{n}\mathrm{Tr}(A_{n}-z{\mathbf{I}}_{n})^{-1}\,,
\]
where $z\in\mathbb{C}^{+}$, and $\mathbf{{I}}_{n}$ is the identity matrix of
order $n$. It is well-known that the Stieltjes transform determines the measure.

The L\'{e}vy distance between two distribution functions $F$ and $G$ is
defined by
\[
L(F,G)=\inf\{\varepsilon>0\ :\ F(x-\varepsilon)-\varepsilon\leq G(x)\leq
F(x+\varepsilon)+\varepsilon\}\,.
\]
We mention that a sequence of distribution functions $F_{n}(x)$ converges to a
distribution function $F(x)$ at all continuity points $x$ of $F$ if and only
if $L(F_{n},G)\rightarrow0.$

Let $(X_{k,\ell})_{(k,\ell)\in{\mathbb{Z}}^{2}}$ be an array of real-valued
random variables, and consider its associated symmetric random matrix
$\mathbf{X}_{n}$ of order $n$ defined by
\begin{equation}
\text{$\big (\mathbf{X}_{n}\big )_{\ell,j}=X_{\ell,j}$ if $1\leq j\leq\ell\leq
n$ and $\big (\mathbf{X}_{n}\big )_{\ell,j}=X_{j,\ell}$ if $1\leq\ell<j\leq
n$}\,\text{.} \label{defbfX}%
\end{equation}
Then, define the symmetric matrix of order $n$ by
\begin{equation}
{\mathbb{X}}_{n}:=\frac{1}{n^{1/2}}\mathbf{X}_{n}. \label{defbbX}%
\end{equation}
The aim of this paper is to study the limiting empirical spectral distribution
function of the symmetric matrix ${\mathbb{X}}_{n}$ defined by \eqref{defbbX}
when the random field $(X_{k,\ell})_{(k,\ell)\in{\mathbb{Z}}^{2}}$ is strictly
stationary given by the following dependence structure: for any $(k,\ell
)\in{\mathbb{Z}}^{2}$,
\begin{equation}
X_{k,\ell}=g(\xi_{k-u,\ell-v}\ ;\ (u,v)\in{\mathbb{Z}}^{2}), \label{defYkl}%
\end{equation}
where $(\xi_{i,j})_{(i,j)\in{\mathbb{Z}}^{2}}$ is an array of i.i.d.
real-valued random variables given on a common probability space
$(\Omega,\mathcal{K},\mathbb{P})$, and $g$ is a measurable function from
${\mathbb{R}}^{{\mathbb{Z}}^{2}}$ to ${\mathbb{R}}$ such that ${\mathbb{E}%
}(X_{0,0})=0$ and $\Vert X_{0,0}\Vert_{2}<\infty$. A representation as in
\eqref{defYkl} includes as special cases, linear as well as many widely used
nonlinear random fields models.

We are interested to establish the weak convergence, on a set of probability
one, of $\nu_{{\mathbb{X}}_{n}}$ to a nonrandom probability measure. This
means that
\begin{equation}
\mathbb{P}(L(F^{{\mathbb{X}}_{n}(\omega)},F)\rightarrow0)=1.
\label{conv in distribution}%
\end{equation}
In the sequel we shall denote this convergence as $F^{{\mathbb{X}}_{n}%
(\omega)}\Rightarrow F$ a.s.

In this paper, for the model defined by (\ref{defbbX}), we shall study the
limit of the type (\ref{conv in distribution})$\ $and specify the limiting
distribution $F(t)$ by giving an equation satisfied by its Stieltjes transform.

Relevant to our result is the notion of spectral density for a weakly
stationary field. In the context of weakly stationary random fields it is
known that, according to Herglotz representation, there exists a unique
measure on $[0,2\pi]^{2}=[0,2\pi]\times\lbrack0,2\pi]$, such that
\[
\gamma_{k,\ell}=\mathrm{cov}(X_{0,0},X_{k,\ell})=\int_{[0,2\pi]^{2}}%
\mathrm{e}^{\mathrm{i}(ku+\ell v)}F(\mathrm{d}u,\mathrm{d}v),\quad\text{for
all}\,\,k,\ell\in\mathbb{Z}\,.
\]
If $F$ is absolutely continuous with respect to the Lebesgue measure
$\lambda\times\lambda$ on $[0,2\pi]^{2}$ then, the Radon-Nikodym derivative
$f$ of $F$ with respect to the Lebesgue measure satisfies%
\[
\gamma_{k,\ell}=\int_{[0,2\pi]^{2}}\mathrm{e}^{\mathrm{i}(ku+\ell
v)}f(u,v)\mathrm{d}u\mathrm{d}v,\quad\text{for all}\,\,k,\ell\in\mathbb{Z}\,.
\]
The function $f(u,v)$ is called \textit{spectral density}.

It should be noted that, by a recent result in Lifshitz and Peligrad \cite{LP}
for random fields defined by \eqref{defYkl}, the spectral density exists. It
is convenient to scale $f(u,v)$ and we define%
\begin{equation}
b(x,y)=(2\pi)^{2}f(2\pi x,2\pi y). \label{defb}%
\end{equation}

One of the main results of this paper is the following theorem, which points
out the relationship between the limiting spectral distribution and the
spectral density.

\begin{theorem}
\label{main}Let $(X_{k,\ell})_{(k,\ell)\in{\mathbb{Z}}^{2}}$ be a real-valued
random field given by \eqref{defYkl} with a spectral density $f(x,y)$. Define
its scaling $b(x,y)$ by (\ref{defb})\ and assume that $\gamma_{k,\ell}%
=\gamma_{\ell,k}$ for all $(k,\ell)\in\mathbb{Z}^{2}$. Then, the convergence
(\ref{conv in distribution}) holds, namely $F^{{\mathbb{X}}_{N}(\omega
)}\Rightarrow F$ a.s., where $F$ is a nonrandom limiting distribution function
whose Stieltjes transform $S(z)$ is uniquely defined by the relations: for
every $z\in{\mathbb{C}}^{+}$
\begin{equation}
S(z)=\int_{0}^{1}g(x,z)\mathrm{d}x,\label{eq0}%
\end{equation}
where, for any $x\in\lbrack0,1],$ $g(x,z)$ is analytic in $z\in{\mathbb{C}%
}^{+}.$ \newline There is $J\subset\lbrack0,1],$ with $\lambda(J)=1$ such
that, for any $x\in J$ and $z\in{\mathbb{C}}^{+},$ $g(x,z)$ satisfies the
equation
\begin{equation}
g(x,z)=-\Big (z+\int_{0}^{1}g(y,z)b(x,y)\mathrm{d}y\Big )^{-1}.\text{
}\label{eq1}%
\end{equation}
Moreover, for any $x\in\lbrack0,1]$ and $z\in{\mathbb{C}}^{+}$
\begin{equation}
\operatorname{Im}g(x,z)>0,\text{ }|g(x,z)|\leq(\operatorname{Im}%
z)^{-1}.\label{eq2}%
\end{equation}

\end{theorem}

This theorem is related to Theorem 3 in \cite{BMP}. The main difference is
that Theorem 3 in \cite{BMP} is obtained under the condition that the
covariances are absolutely summable. This summability condition implies that
the spectral density is continuous and bounded, case known under the name of
short memory. By removing this condition, our Theorem \ref{main} can be
applied to any symmetric random field defined by \eqref{defYkl}, therefore the
memory is not restricted to short memory. The other difference is that
(\ref{eq1}) holds only on a set of Lebesgue $1$ which does not depend on $z$.
But keep in mind that the function given by (\ref{eq1}) is integrated to give
(\ref{eq0}) so $S(z)$ is well determined. Also, this set $J$ can be used to
obtain a version of the spectral density such that equation (\ref{eq1}) holds
for all $x\in\lbrack0,1].$

By using a closely related approach we used to prove Theorem \ref{main}, we
can easily study another symmetrized model based on the random field defined
by \eqref{defYkl}. Instead of $\mathbf{X}_{n}$ defined by (\ref{defbfX}) we
can consider the symmetrized model%
\begin{equation}
\text{$\big (\mathbf{X}_{n}^{\prime}\big )_{k,j}=X_{k,j}+X_{j,k}$ if $1\leq
j,k\leq n$ and }{\mathbb{X}}_{n}^{\prime}:=\frac{1}{n^{1/2}}\mathbf{X}%
_{n}^{\prime}\text{.} \label{deffXprime}%
\end{equation}

For this model, we shall formulate the following result:

\begin{theorem}
\label{cormain}Let $(X_{k,\ell})_{(k,\ell)\in{\mathbb{Z}}^{2}}$ be a
real-valued random field given by \eqref{defYkl} with spectral density
$f(x,y)$. Then, the convergence (\ref{conv in distribution}) holds for
$F^{{\mathbb{X}}_{N}^{\prime}(\omega)}$ to a nonrandom distribution function
$F,$ whose Stieltjes transform, $S(z),$ $z\in{\mathbb{C}}^{+}$, is uniquely
defined by the relations (\ref{eq0}), (\ref{eq1}), (\ref{eq2}), with%
\[
b(x,y)=(2\pi)^{2}(f(2\pi x,2\pi y)+f(2\pi y,2\pi x)).
\]

\end{theorem}

There are certainly connections between the models given by random matrices
(\ref{defbbX}) and (\ref{deffXprime}), as argued in Lemma 19 in \cite{BMP}.
However, Theorem \ref{cormain} does not follow directly from Theorem
\ref{main} since the random field $(X_{k,j}+X_{j,k})_{k,j}$ is no longer
stationary. Their proofs are similar.

It is worth mentioning that a Gaussian random field which has spectral density
is a function of i.i.d., so both Theorem \ref{main} and Theorem \ref{cormain}
apply to this situation.

We shall compare Theorem \ref{cormain} to Theorem 2 in Khorunzhy and Pastur
\cite{KP} (see also in Theorem 17.2.1. in \cite{PS}) concerning Gaussian
random fields. Indeed, if we define%
\[
B(j,k)=\mathrm{cov}(X_{0,0},X_{j,k})+\mathrm{cov}(X_{0,0},X_{k,j}),
\]
then, a straightforward computation shows that
\[
cov(X_{j,k}+X_{k,j},X_{u,v}+X_{v,u})=B(j-u,k-v)+B(j-v,k-u).
\]
Therefore, condition (17.2.3) of Theorem 17.2.1. in \cite{PS} is satisfied.
Note that $B(j,k)=B(k,j)=B(-j,-k).$ As a matter of fact, as noticed in
\cite{BMP}, by a careful analysis of the proof, the condition $B(j,k)=B(j,-k)$
included in (17.2.3) in \cite{PS} can be omitted in the stationary case we
consider. Therefore, for the model treated in Theorem \ref{cormain}, we have
the covariance structure required by Theorem 2 in Khorunzhy and Pastur
\cite{KP}. We can see from these comments that, in the context of the
symmetrized model (\ref{deffXprime}), our Theorem \ref{cormain} extends
Theorem 2 by Khorunzhy and Pastur \cite{KP} (given also in Theorem 17.2.1. in
\cite{PS}) in two directions. The result of Khorunzhy and Pastur \cite{KP}
concerning Gaussian random fields is given under the condition that the
covariances are absolutely summable, implying that the spectral density exists
and is continuous and bounded. We removed this condition in Theorem
\ref{cormain} and for the Gaussian case we can assume only that the spectral
density exists in order to obtain the characterization of the limit. Our
result is also true for random fields which are not necessarily Gaussian, but
are functions of i.i.d.

\begin{remark}
\label{product}If the spectral density of $(X_{k,\ell})_{(k,\ell
)\in{\mathbb{Z}}^{2}}$ has the structure $f(x,y)=u(x)u(y)$ for some real
valued positive function $u,$ then the conclusion of Theorem \ref{main} can be
given in the following simplified form: the convergence
(\ref{conv in distribution}) holds where $F$ is a nonrandom distribution
function whose Stieltjes transform $S(z),$ $z\in\mathbb{C}^{+}$ is given by
the relation
\begin{equation}
S(z)=-\frac{1}{z}(1+v^{2}), \label{star2}%
\end{equation}
where $b(x)=u(2\pi x)$ and $v(z)$ is solution to the equation
\begin{equation}
v(z)=-\int_{0}^{1}\frac{b(y)\mathrm{d}y}{z+b(y)v(z)}\text{,\ }z\in{\mathbb{C}%
}^{+}\,. \label{star1}%
\end{equation}
with $v(z)$ analytic, $\operatorname{Im}v(z)>0,$ and $|v(z)|\leq
(\operatorname{Im}z)^{-1}||X_{0,0}||_{2}.$
\end{remark}

In this form, we can see that one can obtain explicit polynomial equations for
$S(z)$ when $u(x)$ is a positive step function. In particular, if $(X_{k,\ell
})_{(k,\ell)\in{\mathbb{Z}}^{2}}$ is an array of i.i.d. random variables with
zero mean and variance $\sigma^{2}$, then $u$ is constant and $S(z)$ given in
Remark \ref{product} satisfies the equation $\sigma^{2}S^{2}+S+z^{-1}=0$.
Specifying the square root of a complex number as the one with positive
imaginary part, the solution with positive imaginary part is $S=-(z-\sqrt
{z^{2}-4\sigma^{2}})(2\sigma^{2})^{-1}$ which is the well-known Stieltjes
transform of the semicircle law obtained by Wiener \cite{Wiener} (see Lemma
2.11. in \cite{BS}).

In the last section, we are going to provide two examples of random fields,
one linear and one nonlinear, where our results apply. It is remarkable that
for both these examples the only condition required is that the random fields
are well defined in $\mathbb{L}_{2}$.

\section{Proofs}

Before proving the results we shall give some results and facts which will be
used in their proofs.

\subsection{Preliminary considerations}

1. \textbf{Continuity results for Stiletjes transform of limiting spectral
distribution.}

\bigskip

We start by proving a continuity result related to the equation satisfied by
the Fourier transform in the theorems. As mentioned by Khorunzhy and Pastur
\cite{KP}, it should be noted that these equations appear for the first time
in Wegner \cite{W}, in the context of studying $n$-component generalization of
discrete Schr\"{o}dinger equation with a random $n\times n$ matrix. More
precisely, equations (17.2.28) and (17.2.29) in \cite{PS} are comparable to
(3.19) and (3.14) respectively in \cite{W}. Therefore, next proposition has
interest in itself.

\begin{proposition}
\label{cont}Assume that $b_{m}(x,y)$ is a sequence of real positive and
bounded functions, which is convergent in $\mathbb{L}_{1}([0,1]^{2}%
,\mathcal{B}^{2},\lambda^{2})$ to $b(x,y).$ For any $z\in{\mathbb{C}}^{+}$
define
\begin{equation}
S_{m}(z)=\int_{0}^{1}h_{m}(x,z)\mathrm{d}x, \label{eq0m}%
\end{equation}
where $h_{m}(x,z)$ is a solution to the equation
\begin{equation}
h_{m}(x,z)=-\Big (z+\int_{0}^{1}h_{m}(y,z)b_{m}(x,y)\mathrm{d}y\Big )^{-1},
\label{eq1m}%
\end{equation}
with $h_{m}(x,z)$ analytic in $z\in{\mathbb{C}}^{+}$ and for any $x\in
\lbrack0,1],$%
\begin{equation}
\operatorname{Im}h_{m}(x,z)>0,\text{ }|h_{m}(x,z)|\leq(\operatorname{Im}%
z)^{-1}. \label{eq2m}%
\end{equation}
Then, for any $z\in{\mathbb{C}}^{+}$ we have $S_{m}(z)\rightarrow S(z),$ where
$S(z)$ is a Stieltjes transform of a probability measure uniquely determined
by the equations (\ref{eq0}), (\ref{eq1}) and (\ref{eq2}).
\end{proposition}

\textbf{Proof. }It is convenient to represent relation (\ref{eq1m}) in an
equivalent form, by introducing the following transformation. For all
$z\in\mathbb{C}^{+}$ and $x\in\lbrack0,1]$ set
\begin{equation}
h_{m}(x,z)=-(z+\pi_{m}(x,z))^{-1}, \label{defpi}%
\end{equation}
which is well defined since $\operatorname{Im}h_{m}(x,z)>0.$ Equation
(\ref{eq1m}) becomes
\[
\pi_{m}(x,z)=-\int_{0}^{1}\frac{b_{m}(x,s)}{z+\pi_{m}(s,z)}\mathrm{d}s.\,
\]
Also, because%
\[
\pi_{m}(x,z)=\int_{0}^{1}b_{m}(x,s)h_{m}(s,z)\mathrm{d}s
\]
and by the fact that $\operatorname{Im}h_{m}(x,z)>0$ and $b_{m}$ is positive,
it follows that for all $x\in\lbrack0,1]$ and $z\in\mathbb{C}^{+}$ \
\begin{equation}
\operatorname{Im}\pi_{m}(x,z)\geq0. \label{impi}%
\end{equation}
Finally, set $S_{m}(z)=\int_{0}^{1}h_{m}(x,z)\mathrm{d}x.$

Let us show that $\pi_{m}(x,z)$ is pointwise convergent for all $x\in
\lbrack0,1]$ and $z\in\mathbb{C}^{+}.$ Obviously
\begin{gather*}
\pi_{n}(x,z)-\pi_{m}(x,z)=\int_{0}^{1}[\frac{b_{m}(x,s)}{z+\pi_{m}(s,z)}%
-\frac{b_{n}(x,s)}{z+\pi_{n}(s,z)}]\mathrm{d}s=\\
\int_{0}^{1}[\frac{b_{m}(x,s)-b_{n}(x,s)}{z+\pi_{m}(s,z)}+b_{n}(x,s)\frac
{\pi_{n}(s,z)-\pi_{m}(s,z)}{(z+\pi_{m}(s,z))(z+\pi_{n}(s,z))}]\mathrm{d}s.
\end{gather*}
Now, by (\ref{defpi}) and (\ref{eq2m}) it follows that%
\begin{equation}
|z+\pi_{m}(x,z)|=|h_{m}^{-1}(x,z)|\geq\operatorname{Im}z.
\label{boundzipluspi}%
\end{equation}
Also, since $b_{m}$'s are bounded it follows from (\ref{eq1m})\ that we can
find positive constants $K_{m},$ such that
\[
|\pi_{m}(x,z)|\leq\int_{0}^{1}\frac{b_{m}(x,s)}{|z+\pi_{m}(s,z)|}%
\mathrm{d}s\leq\frac{1}{\operatorname{Im}z}K_{m}.
\]
So, for $z\in\mathbb{C}^{+},$ it follows by the above considerations that%
\begin{gather}
\sup_{x\in\lbrack0,1]}|\pi_{n}(x,z)-\pi_{m}(x,z)|\leq\frac{1}%
{\operatorname{Im}z}\int_{0}^{1}|b_{m}(x,s)-b_{n}(x,s)|\mathrm{d}%
s+\label{suppai}\\
\frac{1}{(\operatorname{Im}z)^{2}}\sup_{x\in\lbrack0,1]}|\pi_{n}(x,z)-\pi
_{m}(x,z)|\int_{0}^{1}b_{n}(x,s)\mathrm{d}s<\infty.\nonumber
\end{gather}
By integrating with $x$ on $[0,1]$ we obtain
\begin{gather*}
\sup_{x\in\lbrack0,1]}|\pi_{n}(x,z)-\pi_{m}(x,z)|\leq\frac{1}%
{\operatorname{Im}z}\int_{[0,1]^{2}}|b_{m}(x,s)-b_{n}(x,s)|\mathrm{d}%
s\mathrm{d}x+\\
\frac{1}{(\operatorname{Im}z)^{2}}\sup_{x\in\lbrack0,1]}|\pi_{n}(x,z)-\pi
_{m}(x,z)|\int_{[0,1]^{2}}b_{n}(x,s)\mathrm{d}x\mathrm{d}s.
\end{gather*}
Now, since$\ b_{n}\rightarrow b$ in $\mathbb{L}_{1}[0,1]^{2}$ we have
\begin{equation}
\lim_{n\rightarrow\infty}\int_{[0,1]^{2}}b_{n}(x,s)\mathrm{d}x\mathrm{d}%
s=\int_{[0,1]^{2}}b(x,s)\mathrm{d}x\mathrm{d}s=B. \label{ineqbn}%
\end{equation}
Define the domain $D$ by
\begin{equation}
D=\{z\in\mathbb{C}^{+}:B<\operatorname{Im}z\}. \label{defD}%
\end{equation}
For $z\in D,$ by (\ref{suppai}), (\ref{ineqbn}) and simple algebra, we have
that
\begin{align*}
&  \lim_{m>n\rightarrow\infty}\sup_{x\in\lbrack0,1]}|\pi_{n}(x,z)-\pi
_{m}(x,z)|\\
&  \leq\frac{\operatorname{Im}z}{(\operatorname{Im}z)^{2}-B}\lim
_{m>n\rightarrow\infty}\int_{[0,1]^{2}}|b_{m}(x,s)-b_{n}(x,s)|\mathrm{d}%
s\mathrm{d}x.
\end{align*}
As a consequence, for $z\in D$ we obtain that%
\[
\lim_{m>n\rightarrow\infty}|\pi_{n}(x,z)-\pi_{m}(x,z)|=0\text{ uniformly in
}x,
\]
and by the Lebesgue dominated convergence theorem we also have
\[
\lim_{m>n\rightarrow\infty}\int_{0}^{1}|\pi_{n}(x,z)-\pi_{m}(x,z)|\mathrm{d}%
x=0.
\]
It follows that for $z$ in $D$ and $x\in\lbrack0,1],$ we have both $\pi
_{m}(x,z)$ is convergent pointwise and in $\mathbb{L}_{1}[0,1]$ to a function
we shall denote by
\[
\lim_{m\rightarrow\infty}\pi_{m}(x,z)=:\pi(x,z).
\]

By (\ref{defpi}) we know that $h_{m}(x,z)=-(z+\pi_{m}(x,z))^{-1},$ and by
(\ref{eq2m}), $|h_{m}(x,z)|\leq1/\operatorname{Im}(z).$ Therefore for $z$ in
$D$
\begin{equation}
\lim_{m\rightarrow\infty}\text{ }h_{m}(x,z)=-(z+\pi(x,z))^{-1}=:\text{
}h(x,z), \label{defh1}%
\end{equation}
pointwise and in $\mathbb{L}_{1}[0,1].$ Since for all $x\in\lbrack0,1],$
$(h_{m}(x,z))_{m\geq1}$ are analytic and uniformly bounded functions on all
compacts of $\mathbb{C}^{+}$, by applying Lemma 3 in \cite{GH}, the limit
$h(x,z)$ is also analytic on $D$ for all $x$.

We shall remove now the restriction about $z\in D$ and we shall argue that
actually, for all $x,$ $h_{m}(x,z)$ is convergent on $\mathbb{C}^{+}$ to an
analytic function $g(x,z)$ which coincides with $h(x,z)$ on $D.$ Since
$h_{m}(x,z)$ is a sequence of analytic functions on $\mathbb{C}^{+}$,
uniformly bounded on compacts of $\mathbb{C}^{+}$, by a classical result (see
Lemma 3 in \cite{GH} and references therein) every subsequence of $h_{m}(x,z)$
has a subsequence that converges uniformly on compact subsets of
$\mathbb{C}^{+}$ to an analytic function on $\mathbb{C}^{+}.$ Let us consider
now a convergent subsequence $h_{m^{\prime}}(x,z)\rightarrow g_{1}(x,z)$ with
$g_{1}(x,z)$ analytic on $\mathbb{C}^{+}.$ Clearly on $D$ we have
$g_{1}(x,z)=h(x,z).$ Now if we have another convergent subsequence
$h_{m"}(x,z)\rightarrow g_{2}(x,z)$ with $g_{2}(x,z)$ analytic on
$\mathbb{C}^{+}$, we have $g_{2}(x,z)=h(x,z)$ on $D.$ Therefore, both
$g_{1}(x,z)$ and $g_{2}(x,z)$ are analytic on $\mathbb{C}^{+}$ and
$g_{1}(x,z)-g_{2}(x,z)=0$ on $D.$ It is a well-known fact that two analytic
functions on a connected domain, which coincide on a subset of the domain with
an accumulation point, coincide everywhere on the domain. Therefore we have
that $g_{1}(x,z)=g_{2}(x,z)$ on $\mathbb{C}^{+}$. As a consequence, for each
$x$ fixed, $h_{m}(x,z)\rightarrow g(x,z)$ for all $z\in\mathbb{C}^{+}$ where
$g(x,z)$ is analytic on $\mathbb{C}^{+}$ and coincides with $h(x,z)$ on $D.$

We shall pass now to the limit in the equations (\ref{eq0m}), (\ref{eq1m}) and
(\ref{eq2m}). Since $h_{m}(x,z)$ are bounded, by passing to the limit in
(\ref{eq0m}), we obtain
\[
S(z)=\int_{0}^{1}g(x,z)\mathrm{d}x
\]
for all $z\in\mathbb{C}^{+},$ so (\ref{eq0}) holds$.$

Also, by passing to the limit in (\ref{eq2m}),\ we immediately obtain
\begin{equation}
|g(x,z)|\leq\operatorname{Im}z^{-1}. \label{new}%
\end{equation}
By passing to the limit in equation (\ref{eq1m}) we have%
\begin{equation}
g(x,z)=-\Big (z+\lim_{m\rightarrow\infty}\int_{0}^{1}h_{m}(y,z)b_{m}%
(x,y)\mathrm{d}y\Big )^{-1}. \label{limit}%
\end{equation}
Then, since $\operatorname{Im}h_{m}(x,z)>0$ and $b_{m}(x,y)$ is positive, we
deduce that $\operatorname{Im}g(x,z)>0$ for all $x\in\lbrack0,1]$ and
$z\in\mathbb{C}^{+}.$ It follows that (\ref{eq2})\ holds.

We shall show now that equation (\ref{eq1}) is satisfied for all $x$ in a set
of Lebesgue measure $1$ which does not depend on $z$. With this goal in mind
we shall start from relation (\ref{limit}). Since it gives that for any
$x\in\lbrack0,1]$ and $z\in\mathbb{C}^{+}$ the limit $\lim_{m\rightarrow
\infty}\int_{0}^{1}h_{m}(y,z)b_{m}(x,y)\mathrm{d}y$ exists, in order to verify
(\ref{eq1}) it is enough to find a subsequence ($m^{\prime}$) such that for
all $z\in\mathbb{C}^{+}$ and all $x\in J\subset\lbrack0,1]$ with
$\lambda(J)=1$,%

\begin{equation}
\lim_{m^{\prime}\rightarrow\infty}\int_{0}^{1}h_{m^{\prime}}(y,z)b_{m^{\prime
}}(x,y)\mathrm{d}y=\int_{0}^{1}g(y,z)b(x,y)\mathrm{d}y. \label{limit1}%
\end{equation}
We have
\begin{gather*}
\int_{0}^{1}|h_{m}(y,z)b_{m}(x,y)-g(y,z)b(x,y)|\mathrm{d}y\leq\int_{0}%
^{1}|h_{m}(y,z)-g(y,z)|b(x,y)\mathrm{d}y+\\
\int_{0}^{1}|h_{m}(y,z)||b_{m}(x,y)-b(x,y)|\mathrm{d}y=:I_{m}(x,z)+II_{m}%
(x,z).
\end{gather*}
Note that for all $x\in\lbrack0,1]$ and $z\in\mathbb{C}^{+}$
\[
|h_{m}(y,z)-g(y,z)|b(x,y)\leq\frac{2}{\operatorname{Im}z}b(x,y)
\]
and for $\lambda-$almost all $x\in\lbrack0,1],$ $b(x,y)$ is integrable.
Therefore, by the Lebesgue dominated convergence theorem, $I_{m}%
(x,z)\rightarrow0$ for $x$ in a subset of $[0,1]$ of Lebesgue measure $1$ and
any $z\in\mathbb{C}^{+}.$

By (\ref{eq2m}),%
\[
II_{m}(x,z)\leq\frac{1}{\operatorname{Im}z}\int_{0}^{1}|b_{m}%
(x,y)-b(x,y)|\mathrm{d}y.
\]
Since $b_{m}(x,y)\rightarrow b(x,y)$ in $\mathbb{L}_{1}[0,1]^{2},$ by Fubini
Theorem, $\int_{0}^{1}|b_{m}(x,y)-b(x,y)|\mathrm{d}y$ converges to $0$ in
$\mathbb{L}_{1}[0,1];$ therefore there is a subsequence $(m^{\prime})$ which
does not depend on $z$ such that $\int_{0}^{1}|b_{m^{\prime}}%
(x,y)-b(x,y)|\mathrm{d}y\rightarrow0$ for $\lambda-$almost all $x\in
\lbrack0,1].\ $These considerations show that (\ref{limit1}) holds and
therefore equation (\ref{eq1}) is satisfied for all $z\in\mathbb{C}^{+}$ and
$x\in J\subset\lbrack0,1]$ with $\lambda(J)=1$.

We shall verify now that $S$ is a Stieltjes transform of a probability
measure. Note that since $S$ is a pointwise limit in $\mathbb{C}^{+}$ of
Stieltjes transforms, according to Theorem 1 in Geronimo and Hill \cite{GH},
we have to show that $\ $%
\begin{equation}
\lim_{u\rightarrow\infty}iuS(\mathrm{i}u)=-1. \label{st}%
\end{equation}
A simple computation shows that, by the definition of $g(x,iu),$ we have for
$x\in J$
\[
iug(x,iu)=-\Big (1+\frac{1}{iu}\int_{0}^{1}g(y,iu)b(x,y)\mathrm{d}%
y\Big )^{-1}.
\]
Note that, by (\ref{new}),
\[
|\frac{1}{iu}g(x,iu)|\leq\frac{1}{u^{2}}\rightarrow0\text{ as }u\rightarrow
\infty.
\]
Therefore we can conclude that for all $x\in J$%
\[
iug(x,iu)\rightarrow-1.
\]
Now, again by (\ref{new}) we have $|iug(x,iu)|\leq u/u=1.$ We can apply next
the Lebesgue dominated convergence theorem and obtain%
\[
iuS(iu)=\int_{0}^{1}iyh(x,iu)\mathrm{d}x\,\rightarrow-1.
\]
So, (\ref{st})\ follows, showing that$\ S$ is indeed a Stieltjes transform of
a probability measure with distribution $F$.

It is easy to see that $S(z)$ is uniquely determined by the relations
(\ref{eq0})-(\ref{eq2}). It is convenient to work with the equivalent form of
equation (\ref{eq1}), namely $h(x,z)=-(z+\pi(x,z))^{-1}$. Then, for almost all
$x,$
\begin{equation}
\pi(x,z)=-\int_{0}^{1}\frac{b(x,s)}{z+\pi(s,z)}\mathrm{d}s. \label{pi}%
\end{equation}
If we have two functions $\pi_{1}(x,z)$ and $\pi_{2}(x,z),$ both analytic in
$z,$ we shall write equation (\ref{pi}) for $\pi_{1}(x,z)$ and $\pi_{2}(x,z),$
and by similar manipulations done at the beginning of the proof for $\pi
_{m}(x,z)$ and $\pi_{n}(x,z),$ we get
\[
(1-\frac{B}{(\operatorname{Im}z)^{2}})\sup_{x\in J}|\pi_{1}(x,z)-\pi
_{2}(x,z)|=0\text{ \ for }x\in J.\text{ }%
\]
So, for all $x\in J,$ $\pi_{1}(x,z)=\pi_{2}(x,z)$ and therefore $h_{1}%
(x,z)=-(z+\pi_{1}(x,z))^{-1}=-(z+\pi_{2}(x,z))^{-1}=h_{2}(x,z).$ The
uniqueness follows after we integrate $h_{1}(x,z)$ and $h_{2}(x,z)$ with
respect to $x.$ The proof of Proposition 4 is now complete.

\bigskip

2. \textbf{Facts about universality results for limiting spectral
distribution}

\bigskip

Proposition A below, proved in \cite{BMP}, shows a universality scheme for the
random matrix ${\mathbb{X}}_{n}$ when each $X_{k,\ell}$ is a function of
i.i.d. random variables defined by (\ref{defYkl}). Next, we shall introduce
next a Gaussian random field with the same covariance structure.

Let $(G_{k,\ell})_{(k,\ell)\in{\mathbb{Z}}^{2}}$ be a real-valued centered
Gaussian random field, with covariance function given by
\begin{equation}
{\mathbb{E}}(G_{k,\ell}G_{u,v})={\mathbb{E}}(X_{k,\ell}X_{u,v})\text{ for any
$(k,\ell)$ and $(u,v)$ in ${\mathbb{Z}}^{2}$}\,. \label{egacovfunction}%
\end{equation}
Let $\mathbf{G}_{n}$ be the symmetric random matrix defined by
\begin{equation}
\text{$\big (\mathbf{G}_{n}\big )_{\ell,j}=G_{\ell,j}$ if $1\leq j\leq\ell\leq
n$ and $\big (\mathbf{G}_{n}\big )_{\ell,j}=G_{j,\ell}$ if $1\leq\ell<j\leq
n$}\,\text{.} \label{defbfG}%
\end{equation}
Denote
\begin{equation}
{\mathbb{G}}_{n}=\frac{1}{n^{1/2}}\mathbf{G}_{n}. \label{defbbG}%
\end{equation}

\bigskip

The following is Proposition 1 in\textit{ }\cite{BMP} which shows that the
study of the empirical distribution function of a class of processes which are
functions of i.i.d. random variables can be reduced to the study of a matrix
with Gaussian entries.

\bigskip

\textbf{Proposition A. (}Banna-Merlevede-Peligrad) \textit{Define }%
$(X_{\ell,j})$ by \textit{(\ref{defYkl}), the centered Gaussian random field
}$(G_{k,\ell})$ satisfying (\ref{egacovfunction})$,$\textit{ and the symmetric
matrices }$\mathbb{X}_{n}$\textit{ and }$\mathbb{G}_{n}$\textit{ by
\eqref{defbbX} and \eqref{defbbG} respectively. Then, for any }$z\in C^{+}%
$\textit{, }%
\[
\lim_{n\rightarrow\infty}\big |S_{{\mathbb{X}}_{n}}(z)-{\mathbb{E}%
}\big (S_{{\mathbb{G}}_{n}}(z)\big )\big |=0\ \text{\textit{almost surely}.}%
\]

\bigskip

The following corollary stated in \cite{BMP}, is a direct consequence of
Proposition A together with Theorem B.9 in Bai-Silverstein \cite{BS} :

\bigskip

\textbf{Corollary B}. \textit{Assume that }$\mathbb{X}_{n}$\textit{ and
}$\mathbb{G}_{n}$\textit{ are as in Proposition A. Furthermore, assume that
there exists a nonrandom distribution function }$F$\textit{ such that }%
\begin{equation}
{\mathbb{E}}\big (F^{{\mathbb{G}}_{n}}(t)\big )\rightarrow F(t)\text{
\textit{for all continuity points} }t\in{\mathbb{R}}\text{ \textit{of} }F\,.
\label{exp of F}%
\end{equation}
\textit{Then (\ref{conv in distribution}) holds.}

\bigskip

3. \textbf{Facts about stationary Gaussian fields with spectral density}

\bigskip

Now we mention several facts about stationary Gaussian random fields with
spectral density $f(x,y).$ These facts are also used in \cite{chak} and, for
the case of Gaussian sequences, explained in Ch. 6, Section 6.6. in Varadhan
\cite{Var}.

A centered Gaussian field $(G_{k,\ell})$ has a spectral density $f(x,y)$ if
and only if $(G_{k,\ell})$ is distributed as
\[
G_{k,\ell}=^{d}\sum_{(u,v)\in\mathbb{Z}^{2}}a_{u,v}\xi_{k-u,\ell-v},
\]
with $(\xi_{k,\ell})$ a Gaussian field of i.i.d. random variables centered and
square integrable with variance $1$ and
\begin{equation}
\sum_{(u,v)\in\mathbb{Z}^{2}}a_{u,v}^{2}<\infty. \label{square summable}%
\end{equation}
It is known that
\[
f^{1/2}(x,y)=\frac{1}{2\pi}\left\vert \sum_{(u,v)\in\mathbb{Z}^{2}}%
a_{u,v}\mathrm{e}^{-\mathrm{i}(xu+yv)}\right\vert \ .
\]
Denote
\begin{equation}
G_{k,\ell}^{m}=\sum_{-m\leq u,v\leq m}a_{u,v}\xi_{k-u,\ell-v}. \label{defgm}%
\end{equation}
Let $f_{m}$ be the spectral density of $G_{k,\ell}^{m}$. Then
\[
f_{m}^{1/2}(x,y)=\frac{1}{2\pi}\left\vert \sum_{-m\leq u,v\leq m}%
a_{u,v}\mathrm{e}^{-\mathrm{i}(xu+yv)}\right\vert .
\]
Let us show that $f_{m}^{1/2}(x,y)$ converges in $\mathbb{L}_{2}([0,2\pi
]^{2})$ to $f^{1/2}(x,y).$ Indeed, simple estimates show that%

\begin{align*}
&  \int_{[0,2\pi]^{2}}|f_{m}^{1/2}(x,y)-f^{1/2}(x,y)|^{2}\mathrm{d}%
x\mathrm{d}y\\
&  \leq\frac{1}{(2\pi)^{2}}\int_{[0,2\pi]^{2}}|\sum_{-m\leq u,v\leq m}%
a_{u,v}\mathrm{e}^{-\mathrm{i}(xu+yv)}-\sum_{(u,v)\in\mathbb{Z}^{2}}%
a_{u,v}\mathrm{e}^{-\mathrm{i}(xu+yv)}|^{2}\mathrm{d}x\mathrm{d}y\\
&  =\frac{1}{(2\pi)^{2}}\sum_{(-m\leq u,v\leq m)^{\prime}}a_{u,v}%
^{2}\rightarrow0\text{ as }m\rightarrow\infty.
\end{align*}
This convergence implies that $f_{m}(x,y)$ converges to $f(x,y)$ in
$\mathbb{L}_{1}([0,2\pi]^{2}),$ since, by the Cauchy-Schwarz inequality,%
\begin{gather}
(\int_{[0,2\pi]^{2}}|f_{m}(x,y)-f(x,y)|\mathrm{d}x\mathrm{d}y)^{2}%
\leq\label{convfm}\\
\int_{\lbrack0,2\pi]^{2}}|f_{m}^{1/2}(x,y)-f^{1/2}(x,y)|^{2}\mathrm{d}%
x\mathrm{d}y\int_{[0,2\pi]^{2}}(f_{m}^{1/2}(x,y)+f^{1/2}(x,y))^{2}%
\mathrm{d}x\mathrm{d}y\nonumber\\
\leq2(\mathbb{E}(G_{0,0}^{2})+\mathbb{E(}G_{0,0}^{m})^{2})\int_{[0,2\pi]^{2}%
}(f_{m}^{1/2}(x,y)-f^{1/2}(x,y))^{2}\mathrm{d}x\mathrm{d}y,\nonumber
\end{gather}
and
\[
\lim_{m\rightarrow\infty}\mathbb{E(}G_{0,0}^{m})^{2}=\sum_{u,v}a_{u,v}%
^{2}=\mathbb{E(}G_{0,0}^{2})\text{.}%
\]

\subsection{Proof of Theorem \ref{main}}

The proof has several steps. The idea of proof is that, according to Corollary
\textrm{B,} it is enough to show that (\ref{exp of F}) holds for a Gaussian
random filed with the same spectral density. Then, we shall approximate, as
above, the spectral density $f(x,y)$ of $(G_{k,\ell})$ by the spectral density
$f_{m}(x,y)$ of $(G_{k,\ell}^{m})$ defined by (\ref{defgm}). Since for $m$
fixed the sequence of matrices associated to $(G_{k,\ell}^{m})$ satisfies the
conditions of Theorem 3 in \cite{BMP}, we know how to describe the Stieltjes
transform of the nonrandom limiting empirical spectral distribution $F_{m}$ of
$F^{{\mathbb{G}}_{n}^{m}(\omega)}.$ Then, by using Proposition \ref{cont}, we
show that $F_{m}\Rightarrow F$ as $m\rightarrow\infty,$ where $F$ is a
distribution function of a probability measure, and we shall specify the
equations satisfied by $S$, the Stieltjes transform of $F.$ Finally we show
that $F^{{\mathbb{G}}_{n}(\omega)}\Rightarrow F$ a.s.

\bigskip

\textbf{Step 1.} \textbf{We analyze an associated finite range dependent
random field and a sequence of random matrices.}

Let us construct the random field $(G_{k,\ell}^{m})$ as in (\ref{defgm}) and
consider the random matrix%
\[
\text{$\big (\mathbf{G}_{n}^{m}\big )_{\ell,j}=G_{\ell,j}^{m}$ if $1\leq
j\leq\ell\leq n$ and $\big (\mathbf{G}_{n}^{m}\big )_{\ell,j}=G_{j,\ell}^{m}$
if $1\leq\ell<j\leq n$}\,\text{.}%
\]
Note that, by definition,
\[
\mathrm{cov}(G_{0,0}^{m},G_{k,\ell}^{m})=0\text{ for }k^{2}+\ell^{2}>8m^{2},
\]
and therefore, denoting by $M_{m}=\{(k,\ell)\in{\mathbb{Z}}^{2},$ $k^{2}%
+\ell^{2}\leq8m^{2}\},$%
\begin{equation}
\sum\limits_{(k,\ell)\in{\mathbb{Z}}^{2}}|\mathrm{cov}(G_{0,0}^{m},G_{k,\ell
}^{m})|=\sum\limits_{(k,\ell)\in M_{n}}|\mathrm{cov}(G_{0,0}^{m},G_{k,\ell
}^{m})|<\infty. \label{mdep}%
\end{equation}

By Theorem 3 in \cite{BMP} we conclude that the convergence
(\ref{conv in distribution}), namely $F^{{\mathbb{G}}_{n}^{m}(\omega
)}\Rightarrow F_{m}$ a.s. holds for $F^{{\mathbb{G}}_{n}^{m}(\omega)}$, where
$F_{m}$ is a nonrandom distribution function whose Stieltjes transform
$S_{m}(z)$, $z\in{\mathbb{C}}^{+}$ is uniquely determined by the relations:
(\ref{eq0m})-(\ref{eq2m}), where%
\begin{equation}
b_{m}(x,y)=\sum_{(j,k)\in M_{n}}\mathrm{cov}(G_{0,0}^{m},G_{j,k}%
^{m})\mathrm{e}^{-2\pi\mathrm{i}(xj+yk)}=(2\pi)^{2}f_{m}(2\pi x,2\pi y).
\label{def bn}%
\end{equation}
The last equality above follows, because, (\ref{mdep}) implies that
$f_{m}(x,y)$ is bounded and, by the inversion Fourier formula, the following
representation for the spectral density holds:%
\[
f_{m}(x,y)=\frac{1}{(2\pi)^{2}}\sum_{(j,k)\in M_{n}}\mathrm{cov}(G_{0,0}%
^{m},G_{j,k}^{m})\mathrm{e}^{-\mathrm{i}(xj+yk)}.
\]

\bigskip

\textbf{Step 2. Here we show that }$F_{m}\Rightarrow F$ \textbf{and the
Stieltjes transform of }$F$\textbf{ satisfies the equations of Theorem
\ref{main}}.

Recall that $f_{m}(x,y)$ converges to $f(x,y)$ in $\mathbb{L}_{1}([0,1]^{2}),$
as shown in (\ref{convfm}). Taking into account definition (\ref{defb}), note
that%
\[
B=\int_{[0,1]^{2}}b(x,s)\mathrm{d}x\mathrm{d}s=\int_{[0,2\pi]^{2}%
}f(x,s)\mathrm{d}x\mathrm{d}s=||X_{0,0}||_{2}^{2}.
\]
With the notation $b_{m}(x,y)$ from (\ref{def bn}), note that, after a change
of variables,
\[
\int_{\lbrack0,1]^{2}}|b_{m}(x,y)-b(x,y)|\mathrm{d}x\mathrm{d}y=\int
_{[0,2\pi]^{2}}|f_{m}(x,y)-f(x,y)|\mathrm{d}x\mathrm{d}y
\]
and so $b_{m}(x,y)$ converges in $\mathbb{L}_{1}([0,1]^{2})$ to $b(x,y).$

Now we apply Proposition \ref{cont} and deduce that $S_{m}(z)\rightarrow S(z)$
on $\mathbb{C}^{+}$ and $S(z)$ is a non-random Stieltjes transform of a
probability measure, uniquely determined by the equations (\ref{eq0}),
(\ref{eq1}) and (\ref{eq2}). Furthermore, since $S_{m}(z)\rightarrow S(z)$ on
$\mathbb{C}^{+},$ there is a probability measure with distribution function
$F,$ such that $F_{m}\Rightarrow F$ (see for instance Theorem B.9 in
Bai-Silverstein \cite{BS}).

\bigskip

\textbf{Step 3.} \textbf{Now we show that actually }$F^{{\mathbb{G}}%
_{n}(\omega)}\Rightarrow F$ $a.s.$

This is equivalent to proving that for $z\in C^{+}$ we have $S_{^{{\mathbb{G}%
}_{n}(\omega)}}(z)\rightarrow S(z)$ almost surely. Since $G_{i,j}$ are
functions of $i.i.d.,$ by Proposition A, it is enough to show that for all
$z\in C^{+}$, $\mathbb{E(}S_{^{{\mathbb{G}}_{n}}}(z))\rightarrow S(z).$

By using the triangle inequality,%
\begin{align*}
|\mathbb{E(}S_{^{{\mathbb{G}}_{n}}}(z))-S(z)|  &  \leq|\mathbb{E(}%
S_{^{{\mathbb{G}}_{n}}}(z))-\mathbb{E(}S_{^{{\mathbb{G}}_{n}^{m}}%
}(z))|+|\mathbb{E(}S_{^{{\mathbb{G}}_{n}^{m}}}(z))-S_{m}(z)|\\
&  +|S_{m}(z)-S(z)|.
\end{align*}
By the Step 1 of the proof, $F^{{\mathbb{G}}_{n}^{m}(\omega)}\Rightarrow
F_{m}$ a.s. Therefore, for all $\ z\in C^{+}$, we have $S_{^{{\mathbb{G}}%
_{n}^{m}}}(z)\rightarrow S_{m}(z)$ a.s. But since the Stieltjes transforms are
bounded, we also have $\mathbb{E(}S_{^{{\mathbb{G}}_{n}^{m}}}(z))\rightarrow
S_{m}.$ By the Step 2 of the proof, we know that $S_{m}(z)\rightarrow S(z).$
In order to conclude that $\mathbb{E(}S_{^{{\mathbb{G}}_{n}}}(z))\rightarrow
S(z),$ it suffices to show only that%

\begin{equation}
\lim_{m}\lim\sup_{n}|\mathbb{E(}S_{^{{\mathbb{G}}_{n}}}(z))-\mathbb{E(}%
S_{^{{\mathbb{G}}_{n}^{m}}}(z))|=0\text{ for all}\ z\in C^{+}\text{.}
\label{Bill}%
\end{equation}
By Lemma 2.1 in G\"{o}tze \textit{et al. }\cite{GNT} we have%
\[
|S_{^{{\mathbb{G}}_{n}}}(z)-S_{^{{\mathbb{G}}_{n}^{m}}}(z)|^{2}\leq\frac
{1}{|\operatorname{Im}z|^{4}n}\mathrm{Tr}({\mathbb{G}}_{n}(\omega
)-{\mathbb{G}}_{n}^{m}(\omega))^{2}.
\]
Clearly,%
\begin{gather*}
\mathbb{E}\mathrm{Tr}({\mathbb{G}}_{n}(\omega)-{\mathbb{G}}_{n}^{m}%
(\omega))^{2}\leq\frac{2}{n}\sum_{1\leq j\leq\ell\leq n}\mathbb{E(}G_{\ell
,j}-G_{\ell,j}^{m})^{2}=\\
\leq\frac{2}{n}\sum_{1\leq j\leq\ell\leq n}\mathbb{E(}\sum_{(-m\leq u,v\leq
m)^{\prime}}a_{u,v}\xi_{u-\ell,v-j})^{2}\leq2n\mathbb{(}\sum_{(-m\leq u,v\leq
m)^{\prime}}a_{u,v}^{2}).
\end{gather*}
So, for all $n$ and $m$ we have%
\[
\mathbb{E}(|S_{^{{\mathbb{G}}_{n}}}(z)-S_{^{{\mathbb{G}}_{n}^{m}}}%
(z)|^{2})\leq\frac{2}{|\operatorname{Im}z|^{4}}\sum_{(-m\leq u,v\leq
m)^{\prime}}a_{u,v}^{2}%
\]
and therefore (\ref{Bill}) follows by (\ref{square summable}). The proof of
Theorem \ref{main} is complete.

\bigskip

\textbf{Proof of Theorem \ref{cormain}}\ 

Its proof is a variation of the proof of Theorem \ref{main}. We shall mention
only the differences. As before, we associate to $(X_{k,\ell})_{(k,\ell
)\in{\mathbb{Z}}^{2}}$ a real-valued centered Gaussian random field
$(G_{k,\ell})_{(k,\ell)\in{\mathbb{Z}}^{2}}$ satisfying (\ref{egacovfunction}%
). Let $\mathbf{G}_{n}^{\prime}$ and ${\mathbb{G}}_{n}^{\prime}$ be the
symmetric random matrix defined by
\begin{equation}
\big (\mathbf{G}_{n}^{\prime}\big )_{k,j}=G_{k,j}+G_{j,k}\text{ if }1\leq
j,k\leq n\text{ ; }{\mathbb{G}}_{n}^{\prime}=\frac{1}{n^{1/2}}\mathbf{G}%
_{n}^{\prime}. \label{new g prime}%
\end{equation}
Then, we notice that Proposition A and Corollary B (i.e. Theorem 1 and
Corollary 2 in \cite{BMP}) also hold for ${\mathbb{X}}_{n}^{\prime},$ defined
by (\ref{deffXprime}) and ${\mathbb{G}}_{n}^{\prime},$ defined by
(\ref{new g prime}), replacing ${\mathbb{X}}_{n}$ and ${\mathbb{G}}_{n}.$ To
see this, we have to follow the proof in \cite{BMP} with the following change.
We introduce the following random field: for any integers $k\geq j$ and any
positive integer $\ell$ define\
\[
Y_{k,j}^{\ell}=\mathbb{E}(X_{k,j}|\mathcal{F}_{k,j}^{\ell})+\mathbb{E}%
(X_{j,k}|\mathcal{F}_{j,k}^{\ell}),
\]
where $\mathcal{F}_{k,j}^{\ell}=\sigma(\xi_{u,v};|u-k|\leq\ell,|v-j|\leq
\ell).$ Note that the Euclidian distance $d$ between the two points: $(k,j)$
with $k\geq j$ and $(p,q)$ with $p\geq q,$ is the same as the distance between
the sets $\{(k,j)\cup(j,k)\}$ and $\{(p,q)\cup(q,p)\}.$ It follows that
$Y_{k,j}^{\ell}$ and $Y_{p,q}^{\ell}$ are independent when $d^{2}%
((k,j),(p,q);k\geq j$, $p\geq q)>8\ell^{2}$. By using this remark, all the
arguments in the proof of Proposition A (i.e. Theorem 1 in \cite{BMP}) work unchanged.

Using this new version of Proposition A and Corollary B, as in the proof of
Theorem \ref{main}, we reduce the problem to the study of Gaussian random
matrices selected from stationary Gaussian random fields with the same
spectral density.

We define $G_{\ell,j}^{m}$ by (\ref{defgm}) and%
\[
\big (\mathbf{G}_{n}^{\prime m}\big )_{i,j}=G_{i,j}^{m}+G_{j,i}^{m}\text{ if
}1\leq j,i\leq n\text{ ; }{\mathbb{G}}_{n}^{\prime m}=\frac{1}{n^{1/2}%
}\mathbf{G}_{n}^{\prime m}.
\]
A straightforward computation shows that
\[
cov(G_{j,k}^{m}+G_{k,j}^{m},G_{u,v}^{m}+G_{v,u}^{m})=B_{m}(j-u,k-v)+B_{m}%
(j-v,k-u),
\]
where%
\[
B_{m}(j,k)=\mathrm{cov}(G_{0,0}^{m},G_{j,k}^{m})+\mathrm{cov}(G_{0,0}%
^{m},G_{k,j}^{m}).
\]
Therefore condition (17.2.3) of Theorem 17.2.1. in \cite{PS} is satisfied. By
Theorem 17.2.1. in \cite{PS} we conclude that the convergence $F^{{\mathbb{G}%
}_{n}^{\prime m}(\omega)}\Rightarrow F_{m}^{\prime}$ a.s. holds, where
$F_{m}^{\prime}$ is a nonrandom distribution function whose Stieltjes
transform $S_{m}(z)$, $z\in{\mathbb{C}}^{+}$ is uniquely defined for
$z\in\mathbb{C}^{+}$ by the relations (\ref{eq0m})-(\ref{eq2m}) with%
\[
b_{m}(x,y)=(2\pi)^{2}(f_{m}(2\pi x,2\pi y)+f_{m}(2\pi y,2\pi x)).
\]
From here on, the proof of Theorem \ref{cormain}\ is identical to the proof of
Theorem \ref{main}.

\bigskip

\textbf{Proof of Remark \ref{product}. }

By the conditions of this remark, note that $b(x,y)=t(x)t(y)$, with
$t(x)=(2\pi u)(2\pi x),$ which together with (\ref{eq1}) gives%
\begin{equation}
g(x,z)=-(z+t(x)v(z))^{-1},\text{ } \label{star}%
\end{equation}
where
\[
v(z)=\int_{0}^{1}g(y,z)t(y)\mathrm{d}y.
\]
By multiplying (\ref{star}) by $t(x)$ and integrating with $x$ we get the
equation (\ref{star1}).

Now from equation (\ref{eq0}) we have
\[
S(z)=-\int_{0}^{1}\frac{\mathrm{d}x}{z+t(x)v(z)}.
\]
Multiplying equation (\ref{star1}) by $v(z)$ and adding it to $zS(z)$\ we
obtain $v^{2}(z)+zS(z)=-1$, leading to (\ref{star2}). Next, by (\ref{eq2}), we
notice that $v(z)$ is analytic, has strictly positive imaginary part and
\[
|v(z)|\leq\frac{1}{\operatorname{Im}z}\int_{0}^{1}t(y)\mathrm{d}y=\frac
{1}{\operatorname{Im}z}\int_{0}^{2\pi}u(x)\mathrm{d}x=\frac{1}%
{\operatorname{Im}z}||X_{0,0}||_{2}.
\]

\section{Examples}

\subsection{Linear processes}

Let $(a_{k,\ell})_{(k,\ell)\in{\mathbb{Z}}^{2}}$ be a double indexed sequence
of real numbers such that
\begin{equation}
\sum_{(k,\ell)\in{\mathbb{Z}}^{2}}a_{k,\ell}^{2}<\infty\,, \label{longmem}%
\end{equation}
and let $(X_{u,v})_{(u,v)\in{\mathbb{Z}}^{2}}$ be the linear random field
defined by: for any $(u,v)\in{\mathbb{Z}}^{2}$,
\begin{equation}
X_{u,v}=\sum_{(k,\ell)\in{\mathbb{Z}}^{2}}a_{k,\ell}\xi_{k+u,\ell+v}\,,
\label{linproc}%
\end{equation}
where the variables $(\xi_{k,\ell})$ are i.i.d. centered and square
integrable. Note that $X_{u,v}$ is well defined in $\mathbb{L}_{2}$ if and
only if (\ref{longmem}) holds. The spectral density is defined by
\[
f(x,y)=\frac{1}{(2\pi)^{2}}\left\vert \sum_{(u,v)\in\mathbb{Z}^{2}}%
a_{u,v}\mathrm{e}^{-\mathrm{i}(xu+yv)}\right\vert ^{2}.
\]
We can apply our results to the random matrices associated to the linear
random field and obtain the following corollary, describing the limiting
spectral distribution.

\begin{corollary}
Assume that condition (\ref{longmem}) is satisfied and $(X_{u,v}%
)_{(v,v)\in{\mathbb{Z}}^{2}}$ is defined by (\ref{linproc}). Then $X_{u,v}$'s
are well defined in $\mathbb{L}_{2}$ and the conclusion of Theorem
\ref{cormain} holds. If $a_{k,\ell}=a_{\ell,k}$ for any $(k,\ell
)\in{\mathbb{Z}}^{2},$ then the conclusion of Theorem \ref{main}\ holds. If
$a_{k,\ell}=a_{k}a_{\ell}$ for some real numbers $a_{k}$ with $\sum
_{k\in{\mathbb{Z}}}a_{k}^{2}<\infty$, then Remark \ref{product} applies.
\end{corollary}

\subsection{Volterra-type processes}

Other classes of stationary random fields having the representation
\eqref{defYkl} are Volterra-type processes which play an important role in the
nonlinear system theory. For any ${\mathbf{k}}=(k_{1},k_{2})\in{\mathbb{Z}%
}^{2}$, define the Volterra-type expansion as follows:
\begin{equation}
X_{{\mathbf{k}}}=\sum_{{\mathbf{u}},{\mathbf{v}}\in{\mathbb{Z}}^{2}%
}b_{{\mathbf{u}},{\mathbf{v}}}\xi_{{\mathbf{k}}-{\mathbf{u}}}\xi_{{\mathbf{k}%
}-{\mathbf{v}}}\,, \label{Volterra}%
\end{equation}
where $b_{{\mathbf{u}},{\mathbf{v}}}$ are real numbers satisfying
\begin{equation}
b_{{\mathbf{u}},{\mathbf{v}}}=0\text{ if }{\mathbf{u}}={\mathbf{v}},\text{
}\sum_{{\mathbf{u}},{\mathbf{v}}\in{\mathbb{Z}}^{2}}b_{{\mathbf{u}%
},{\mathbf{v}}}^{2}<\infty\,, \label{condcoeffvolt}%
\end{equation}
and $(\xi_{{\mathbf{k}}})_{\mathbf{k}\in{\mathbb{Z}}^{2}}$ is an i.i.d. random
field of centered and square integrable variables. Under the above conditions,
the random field $X_{{\mathbf{k}}}$ exists, is stationary, centered and square
integrable. By \cite{LP}, the field has spectral density since it is a
function of i.i.d. The covariance structure is given by: for any ${\mathbf{u}%
}=(u_{1},u_{2})$ and ${\mathbf{v}}=(v_{1},v_{2})$ in ${\mathbb{Z}}^{2}$,
\begin{equation}
\gamma_{{\mathbf{k}}}=\Vert\xi_{0,0}\Vert_{2}^{4}\sum_{{\mathbf{u}%
},{\mathbf{v}}\in{\mathbb{Z}}^{2}}b_{{\mathbf{u}},{\mathbf{v}}}(b_{{\mathbf{u}%
}+{\mathbf{k}},{\mathbf{v}}+{\mathbf{k}}}+b_{{\mathbf{v}}+{\mathbf{k}%
},{\mathbf{u}}+{\mathbf{k}}})\ \text{ for any ${\mathbf{k}}\in{\mathbb{Z}}%
^{2}$}. \label{defgammakvolt}%
\end{equation}

Therefore the following corollary holds:

\begin{corollary}
\label{propVolt1} Assume that condition (\ref{condcoeffvolt}) is satisfied and
$(X_{{\mathbf{k}}})_{{\mathbf{k}}\in{\mathbb{Z}}^{2}}$ is defined by
(\ref{Volterra}). Then the conclusion of Theorem \ref{cormain} holds for the
model given by (\ref{deffXprime}).
\end{corollary}

If we impose additional symmetry conditions to the coefficients
$b_{{\mathbf{u}},{\mathbf{v}}}$ defining the Volterra random field
\eqref{Volterra}, we can derive the limiting spectral distribution of its
associated symmetric matrix ${\mathbb{X}}_{n}$ defined by \eqref{defbbX}.

If for any ${\mathbf{u}}=(u_{1},u_{2})$ and ${\mathbf{v}}=(v_{1},v_{2})$ in
${\mathbb{Z}}^{2}$,$\,b_{{\mathbf{u}},{\mathbf{v}}}=b_{u_{_{1}},v_{_{1}}%
}b_{u_{_{2}},v_{_{2}}}$ the conclusion of Theorem \ref{main} holds.

\bigskip

\textbf{Acknowledgement.} The research for this paper was supported in part by
a Charles Phelps Taft Memorial Fund grant, and the NSF grant DMS-1512936. The
authors would like to thank the anonymous reviewer for valuable comments,
which improved the presentation of the paper.


\begin{thebibliography}{99}                                                                                               %


\bibitem {AZ}Anderson, G. and Zeitouni, O. (2008). A law of large numbers for
finite-range dependent random matrices. \textit{Comm. Pure Appl. Math.}
\textbf{61} 1118-1154.

\bibitem {BS}Bai, Z. and Silverstein, J.W. (2010). \textit{Spectral analysis
of large dimensional random matrices}. Springer, New York, second edition.

\bibitem {BZ}Bai, Z. and Zhou, W. (2008). Large sample covariance matrices
without independence structures in columns. \textit{Statist. Sinica}
\textbf{18} 425-442.

\bibitem {BMP}Banna, M., Merlev\`{e}de, F. and Peligrad, M. (2015). On the
limiting spectral distribution for a large class of random matrices with
correlated entries. \textit{Stoc. Proc. Appl.} \textbf{125 }2700-2726.

\bibitem {BK}Boutet de Monvel, A. and Khorunzhy, A. (1999). On the Norm and
Eigenvalue Distribution of Large Random Matrices. \textit{Ann. Probab.}
\textbf{27} 913-944.

\bibitem {BKV}Boutet de Monvel, A. Khorunzhy, A. and Vasilchuk, V. (1996).
Limiting eigenvalue distribution of random matrices with correlated entries.
\textit{Markov Process. Related Fields} \textbf{2} 607-636.

\bibitem {chak}Chakrabarty A., Hazra R.S. and Sarkar D. (2014). From random
matrices to long range dependence. \textit{arXiv:1401.0780}.

\bibitem {GH}Geronimo, J. and Hill, T. (2003). Necessary and sufficient
condition that the limit of Stieltjes transforms is a Stieltjes transform.
\textit{J. Approx. Theory} \textbf{121} 54--60.

\bibitem {GNT}G\"{o}tze, F., Naumov, A. and Tikhomirov A. (2013). Semicircle
law for a class of random matrices with dependent entries.
\textit{arXiv:1211.0389v2}.

\bibitem {HLN}Hachem, W., Loubaton, P. and J. Najim (2005). The empirical
eigenvalue distribution of a Gram matrix: from independence to stationarity,
\textit{Markov Process. Related Fields} \textbf{11} 629--648.

\bibitem {KP}Khorunzhy, A. and Pastur, L. (1994). On the eigenvalue
distribution of the deformed Wigner ensemble of random matrices. In: V. A.
Marchenko (ed.), \textit{Spectral Operator Theory and Related Topics}, Adv.
Soviet Math. \textbf{19}, Amer. Math. Soc., Providence, RI, 97-127.

\bibitem {K}Khorunzhy, A. (1996). Eigenvalue distribution of large random
matrices with correlated entries.\textit{ Matematich. Fizika, Analiz i
Geometria} \textbf{3} 80-101.

\bibitem {LP}Lipshitz, M. and Peligrad, M. (2015). On the spectral density of
stationary processes and random fileds. Zapiski Nauchnyh Seminarov POMI, vol.
441. Probability and Statistics 22 (editors A.N.Borodin, M.A.Lifshits,
A.Yu.Zaitsev) 274-285.

\bibitem {MP}Merlev\`{e}de, F. and Peligrad, M. (2014). On the empirical
spectral distribution for matrices with long memory and independent rows.
\textit{arXiv: 1406.1216}

\bibitem {PS}Pastur, L. and Shcherbina, M. (2011). \textit{Eigenvalue
distribution of large random matrices}. Mathematical Surveys and Monographs,
\textbf{171}. American Mathematical Society, Providence, RI.

\bibitem {Var}Varadhan, S.R.S. (2001). \textit{Probability theory}, Courant
lecture notes \textbf{7.} American Mathematical Society.

\bibitem {W}Wegner, F. (1979). Disordered system with n orbitals per site:
n$=\infty$ limit. \textit{Phys. Rev. B} \textbf{19} 783--792.

\bibitem {Wiener}Wigner E.P. (1958). On the distribution of the roots of
certain symmetric matrices.\textit{ Ann. of Math.} \textbf{67} 325-327.

\bibitem {Yao}Yao, J. (2012). A note on a Mar\u{c}enko-Pastur type theorem for
time series. \textit{Statist. Probab. Lett.} \textbf{82} 22-28.
\end{thebibliography}
\end{document}